# Existence of n-cycles and border-collision bifurcations in piecewise-linear continuous maps with applications to recurrent neural networks

Z. Monfared* · D. Durstewitz

**Abstract** Piecewise linear recurrent neural networks (PLRNNs) form the basis of many successful machine learning applications for time series prediction and dynamical systems identification, but rigorous mathematical analysis of their dynamics and properties is lagging behind. Here we contribute to this topic by investigating the existence of n-cycles ($n \geq 3$) and border-collision bifurcations in a class of n-dimensional piecewise linear continuous maps which have the general form of a PLRNN. This is particularly important as for one-dimensional maps the existence of 3-cycles implies chaos. It is shown that these n-cycles collide with the switching boundary in a border-collision bifurcation, and parametric regions for the existence of both stable and unstable n-cycles and border-collision bifurcations will be derived theoretically. We then discuss how our results can be extended and applied to PLRNNs. Finally, numerical simulations demonstrate the implementation of our results and are found to be in good agreement with the theoretical derivations. Our findings thus provide a basis for understanding periodic behavior in PLRNNs, how it emerges in bifurcations, and how it may lead into chaos.

## 1 Introduction

A piecewise smooth discrete-time dynamical system is a discrete-time map whose state space is split into two or more components (sub-regions) by some discontinuity borders or switching manifolds, such that in each sub-region there is a different functional form of the map [2,3,7,26]. Piecewise smooth (PWS) maps have received growing attention in recent years, as they have a wide range of applications in various areas such as neural dynamics, switching circuits, or impacting mechanical systems [23]. One important type of PWS map is a piecewise linear continuous map, which is continuous but has some discontinuities in its Jacobian matrix across the switching boundaries. Piecewise linear recurrent neural networks (PLRNNs), which build on so-called 'Rectified Linear Units (ReLU)', $\phi(z) = \max(z, 0)$, as the network's nonlinear activation function, are one example

Z. Monfared
Department of Theoretical Neuroscience, Central Institute of Mental Health, Medical Faculty Mannheim, University of Heidelberg, Mannheim, Germany.
E-mail: zahra.monfared@zi-mannheim.de

D. Durstewitz[1,2]
[1]Department of Theoretical Neuroscience, Central Institute of Mental Health, Medical Faculty Mannheim, Heidelberg University, Mannheim, Germany.
[2]Faculty of Physics and Astronomy, Heidelberg University, Heidelberg, Germany.
E-mail: daniel.durstewitz@zi-mannheim.de

∗ Corresponding author

of such maps. In general, RNNs are the standard these days in machine learning for processing sequential, time-series information, due to their success in domains rich in temporal structure like natural language processing [15,27], prediction of consumer behavior [16], movement trajectories [19], or identification of dynamical systems from experimental data [14]. ReLU-based RNNs are particularly popular as they allow for highly efficient inference and training algorithms that exploit their piecewise linear structure [17,21,6,14]. To understand the representational and computational capabilities of these systems in the various application areas indicated above, it is important to study them more systematically from a mathematical, dynamical systems perspective. In particular, periodic motion forms an integral part of many natural data domains (like speech or movement signals), and periodic orbits of various orders have been extracted from brain signals recently using PLRNNs [14].

There are different types of bifurcations in PWS dynamical systems, notably bifurcations that occur because of the existence of the discontinuity boundaries. These form the class of discontinuity-induced bifurcations that only exist in PWS systems [4]. In particular, border-collision bifurcations, or C-bifurcations, have many applications in engineering, computational neuroscience, biology, economics and the social sciences, [11]. They arise when either fixed points or periodic points of a PWS map collide with one of the switching boundaries at a critical value of the bifurcation parameter [3, 4]. Feigin [8–10] was among the first to study analytical conditions for border-collision bifurcations of fixed points and period-2 orbits in piecewise linear (PWL) continuous maps. In [3] some theoretical results on the existence of period-2 orbits in n-dimensional PWS continuous maps with one bifurcation parameter were presented. Subsequently, [13] provided a description of border-collision bifurcations in one-dimensional discontinuous maps, as well as some conditions for the creation and stability of different periodic motions and chaos. Later, Higham et al. [12] considered the occurrence of period-1 and period-2 fixed points in n-dimensional PWL maps with a gap. In [7] the existence of period-1 and period-2 orbits during a border collision bifurcation was discussed for n-dimensional PWL discontinuous maps with two parameters. Attractors and border collision bifurcations of the skew tent map, which as will be shown below is tightly related to the system studied here under certain conditions, have been studied in depth in [1,25]. Very recently, Patra [22] investigated the coexistence of a period-2 orbit, a period-3 orbit, and an unstable chaotic orbit for some parameter values of a 3D PWL normal form map (see also [14] for similar numerical observations in PLRNNs inferred from data).

However, since high-dimensional PWL maps, of which PLRNNs are an example, with an exponentially growing number of discontinuity boundaries are hard to handle, here we consider such systems locally in the neighborhood of only one switching manifold. Then, for an $n$-dimensional system, there is one $s$ in $\{1, \ldots, n\}$ such that the system is defined by

$$Z^{(k+1)} = F_\mu(Z^{(k)}) = \begin{cases} A_1 Z^{(k)} + \mu\, h; & z_s^{(k)} \leq 0 \\ A_2 Z^{(k)} + \mu\, h; & z_s^{(k)} \geq 0 \end{cases}. \tag{1}$$

Although (1) is a reduced system, by adding some assumptions we can extend the obtained results to PLRNNs more generally. Up to now, conditions for the existence of period-$k$ orbits ($k \geq 3$) have been established only in one- and two-dimensional PWS maps, but not for $n$-dimensional PWS maps more generally [7,22,25]. In fact, the existence of n-cycles ($n \geq 3$) and chaotic orbits for n-dimensional maps has been an open problem so far [7]. Here we work out the conditions for the existence, related stability and border collision bifurcation of n-cycles ($n \geq 3$), with symbolic sequence $RL^{n-1}$, for a class of n-dimensional PWL maps with one discontinuity boundary, i.e. system (1) with a special form of matrices $A_1$ and $A_2$. We will find parametric regions for the occurrence of stable and unstable n-cycles in such maps. Our theoretical results reveal that these periodic orbits lie precisely on the switching border for a specific value of some system parameters independent of $\mu$, implying that the system undergoes a border-collision bifurcation. We note that for one-dimensional PWL continuous maps the existence of 3-cycles implies the existence of chaos [18]. Finally, we show that our findings can be applied to PLRNNs.

This paper is organized as follows: Section 2 provides some mathematical preliminaries. In section 3, a special form for the matrices $A_1$ and $A_2$ of system (1) is obtained for which the existence, stability and border collision bifurcations of all the n-cycles $RL^{n-1}(n \geq 3)$ will be examined. Section 4 discusses the application of our obtained results to PLRNNs. In section 5, finally, numerical simulations are provided that indicate how to apply our findings.

## 2 Preliminaries

Consider the PWL map (1) on $\mathbb{R}^n$, where $z_s^{(k)}$ ($1 \leq s \leq n$) is the $s$-th component of $Z^{(k)} = (z_1^{(k)}, z_2^{(k)}, ..., z_n^{(k)})^T \in \mathbb{R}^n$, $A_1 = [a_{ij}^{(1)}]$ and $A_2 = [a_{ij}^{(2)}], i,j = 1, 2, \cdots, n$, are $n \times n$ matrices with real entries, $h = (h_1, h_2, ..., h_n)^T \in \mathbb{R}^n$, and $\mu \in \mathbb{R}$ is a real valued parameter of the system. Assume that matrices $A_1$ and $A_2$ are identical except for the s-th column, i.e. $a_{ij}^{(1)} = a_{ij}^{(2)}$ if $j \neq s$. Furthermore, let us denote the discontinuity boundary of map (1) by $\Sigma$, and the two subregions separated through this boundary by $S^-$ and $S^+$:

$$S^- = \{Z^{(k)} = (z_1^{(k)}, z_2^{(k)}, ..., z_n^{(k)})^T \in \mathbb{R}^n;\ H(Z^{(k)}) = z_s^{(k)} < 0\}, \quad (2)$$

$$\Sigma = \{Z^{(k)} = (z_1^{(k)}, z_2^{(k)}, ..., z_n^{(k)})^T \in \mathbb{R}^n;\ H(Z^{(k)}) = z_s^{(k)} = 0\}, \quad (3)$$

$$S^+ = \{Z^{(k)} = (z_1^{(k)}, z_2^{(k)}, ..., z_n^{(k)})^T \in \mathbb{R}^n;\ H(Z^{(k)}) = z_s^{(k)} > 0\}, \quad (4)$$

where the scalar function $H : \mathbb{R}^n \to \mathbb{R}$, with $H(Z^{(k)}) = z_s^{(k)}$ has nonvanishing gradient. Then, we can rewrite map (1) as

$$Z^{(k+1)} = F_\mu(Z^{(k)}) = \begin{cases} F_\mu^-(Z^{(k)}) = A_1 Z^{(k)} + \mu\, h; & Z^{(k)} \in \bar{S}^- \\ F_\mu^+(Z^{(k)}) = A_2 Z^{(k)} + \mu\, h; & Z^{(k)} \in \bar{S}^+ \end{cases}. \quad (5)$$

**Theorem 1 (Period three implies chaos)** *Suppose that $F : I \to I$ is a continuous map with $I \subset \mathbb{R}$. If $F$ has a period-3 orbit, then $F$ is chaotic.*

*Proof* See [18].

## 3 Periodic orbits and bifurcations

Consider PWL map (5) and, without loss of generality, assume that $s = 1$. Then, $A_1$ and $A_2$ are two $n \times n$ matrices that only differ in their first column. All results of this section can be proven analogously for $s \neq 1$.

**Lemma 1** *Let $A_1$ and $A_2$ be two $n \times n$ matrices that differ only in their first column, i.e. $a_{ij}^{(1)} = a_{ij}^{(2)}$ if $j \neq 1$. Then, there are special forms for the matrices $A_1$ and $A_2$ for which $A_1^k$ and $A_2^k$ also differ only in their first column, for all $k = 1, 2, 3, \cdots$. In this case, $A_1^k$, $A_2^k$ and $A_1^{k_1} A_2^{k_2} A_1^{k_3} A_2^{k_4} \ldots A_1^{k_{n-1}} A_2^{k_n} (k_1 + k_2 + k_3 + k_4 + \ldots + k_n = k)$, are also equal except for the first column.*

*Proof* Suppose that $A_1$ and $A_2$ are two $n \times n$ matrices that differ in their first column. Therefore, they can be partitioned in the following way:

$$A_1 = \left(\begin{array}{c|c} a & \vec{c}^T \\ \hline \vec{b} & A \end{array}\right), \qquad A_2 = \left(\begin{array}{c|c} d & \vec{c}^T \\ \hline \vec{e} & A \end{array}\right), \quad (6)$$

such that $A \in \mathbb{R}^{(n-1)\times(n-1)}$ and $\vec{c}, \vec{b}, \vec{e} \in \mathbb{R}^{n-1}$. In this case, $A_1^2$ and $A_2^2$ can be written as

$$A_1^2 = \left(\begin{array}{c|c} a^2 + \vec{c}^T \vec{b} & a\vec{c}^T + \vec{c}^T A \\ \hline a\vec{b} + A\vec{b} & \vec{b}\vec{c}^T + A^2 \end{array}\right), \qquad A_2^2 = \left(\begin{array}{c|c} d^2 + \vec{c}^T \vec{e} & d\vec{c}^T + \vec{c}^T A \\ \hline d\vec{e} + A\vec{e} & \vec{e}\vec{c}^T + A^2 \end{array}\right). \qquad (7)$$

Thus, $A_1^2$ and $A_2^2$ also differ only in their first columns iff

$$\begin{cases} a\vec{c}^T + \vec{c}^T A = d\vec{c}^T + \vec{c}^T A \\ \vec{b}\vec{c}^T + A^2 = \vec{e}\vec{c}^T + A^2 \end{cases} \implies \begin{cases} (a-d)\vec{c}^T = 0 \\ (\vec{b} - \vec{e})\vec{c}^T = 0 \end{cases}. \qquad (8)$$

Since the first columns of $A_1$ and $A_2$ are different, the last assertion holds if $\vec{c} = \vec{0}$. In this case we have:

$$A_1 = \left(\begin{array}{c|c} a & \vec{0}^T \\ \hline \vec{b} & A \end{array}\right), \qquad A_2 = \left(\begin{array}{c|c} d & \vec{0}^T \\ \hline \vec{e} & A \end{array}\right), \qquad (9)$$

$$A_1^2 = \left(\begin{array}{c|c} a^2 & \vec{0}^T \\ \hline a\vec{b} + A\vec{b} & A^2 \end{array}\right), \qquad A_2^2 = \left(\begin{array}{c|c} d^2 & \vec{0}^T \\ \hline d\vec{e} + A\vec{e} & A^2 \end{array}\right), \qquad (10)$$

$$A_1^3 = \left(\begin{array}{c|c} a^3 & \vec{0}^T \\ \hline * & A^3 \end{array}\right), \qquad A_2^3 = \left(\begin{array}{c|c} d^3 & \vec{0}^T \\ \hline * & A^3 \end{array}\right), \qquad (11)$$

$$\vdots \qquad \qquad \vdots$$

$$A_1^k = \left(\begin{array}{c|c} a^k & \vec{0}^T \\ \hline * & A^k \end{array}\right), \qquad A_2^k = \left(\begin{array}{c|c} d^k & \vec{0}^T \\ \hline * & A^k \end{array}\right). \qquad (12)$$

This means that $A_1^k$ and $A_2^k$ ($k = 1, 2, 3, \cdots$) differ only in their first columns. Likewise, in this case $A_1^k$, $A_2^k$ and $A_1^{k_1} A_2^{k_2} A_1^{k_3} A_2^{k_4} \ldots A_1^{k_{n-1}} A_2^{k_n} (k_1 + k_2 + k_3 + k_4 + \ldots + k_n = k)$ are also equal except for the first column, which completes the proof.

Now let us denote $Z = (x, Y)$, where $x \in \mathbb{R}$ and $Y \in \mathbb{R}^{n-1}$ collects the $(n-1)$ remaining elements. Also assume that

$$\mu h = (\mu h_1, \mu h_2, \cdots, \mu h_n)^T := (\hat{\mu}, h_Y)^T, \qquad (13)$$

where $\hat{\mu} := \mu h_1$, and $h_Y$ contains all the other vector components. It will be shown that only the term $\hat{\mu}$ is crucial. Then, system (5) with $A_1$ and $A_2$ given by (9) can be written as

$$Z^{(k+1)} = F(Z^{(k)}) = \begin{cases} F^-(Z^{(k)}); & x^{(k)} \leq 0 \\ F^+(Z^{(k)}); & x^{(k)} \geq 0 \end{cases}, \qquad (14)$$

where

$$F^-(Z^{(k)}) = \begin{cases} a\, x^{(k)} + \hat{\mu} \\ \vec{b}\, x^{(k)} + AY^{(k)} + h_Y \end{cases}, \qquad F^+(Z^{(k)}) = \begin{cases} d\, x^{(k)} + \hat{\mu} \\ \vec{e}\, x^{(k)} + AY^{(k)} + h_Y \end{cases}. \qquad (15)$$

We remark that the first line in (15) yields the skew tent map ([1, 25]; see also below), while the extended multi-dimensional system is introduced to embed this map within a PLRNN ([14, 24]).

3.1 Existence and border collision bifurcations of n-cycles, $n \geq 3$

Here, using an idea of one of the referees, rewriting system (14), first we will investigate the existence of 3-cycles, with symbolic sequence $RL^2$ (and their complementary cycles with symbolic sequence $RLR$), as well as border collision bifurcations of these cycles. Afterwards, the existence, stability and border collision bifurcations of all n-cycles $RL^{n-1}$ (and $RL^{n-2}R$), $n > 3$ will be studied.

Consider the continuous PWL system in (14)-(15). We can rewrite it in the form

$$x^{(k+1)} = f(x^{(k)}) = \begin{cases} f^-(x^{(k)}) = ax^{(k)} + \hat{\mu}; & x^{(k)} \leq 0 \\ f^+(x^{(k)}) = dx^{(k)} + \hat{\mu}; & x^{(k)} \geq 0 \end{cases}, \quad \text{(Skew tent map)} \quad (16)$$

$$Y^{(k+1)} = G(x^{(k)}, Y^{(k)}) = \begin{cases} G^-(x^{(k)}, Y^{(k)}) = \overrightarrow{b} x^{(k)} + AY^{(k)} + h_Y; & x^{(k)} \leq 0 \\ G^+(x^{(k)}, Y^{(k)}) = \overrightarrow{e} x^{(k)} + AY^{(k)} + h_Y; & x^{(k)} \geq 0 \end{cases}, \quad (17)$$

where the one-dimensional system (16) is the skew tent map whose dynamics and bifurcations are well studied (see for example [1,25]). Indeed, for the system (16)-(17), it will be shown that the existence of the cycles depend on this 1D map (16).

Let $\hat{\mu} > 0$ and consider a 3-cycle $RL^2$ with periodic points $Z_1, Z_2$ and $Z_3$. Clearly, a necessary condition for the existence of a 3-cycle of (14) is that the related $x$-components, $x_1, x_2$ and $x_3$, form a cycle of the skew tent map (16). Therefore, $f^- \circ f^- \circ f^+(x_1) = x_1$, and from (16) we have

$$x_1 = \frac{1+a+a^2}{1-a^2d} \hat{\mu}, \quad a^2d \neq 1, \quad (18)$$

such that $x_1 > 0$. Moreover

$$x_2 = f^+(x_1) = \frac{1+d+ad}{1-a^2d} < 0, \quad x_3 = f^-(x_2) = \frac{1+a+ad}{1-a^2d} < 0. \quad (19)$$

On the other hand, to obtain the vector $Z_1 = (x_1, Y_1)$, it is sufficient to solve system (17) for $Y_1$. To this end, from

$$Y_2 = \overrightarrow{e} x_1 + AY_1 + h_Y,$$
$$Y_3 = \overrightarrow{b} x_2 + AY_2 + h_Y = \overrightarrow{b} x_2 + x_1 A \overrightarrow{e} + A^2 Y_1 + Ah_Y + h_Y,$$
$$Y_1 = \overrightarrow{b} x_3 + AY_3 + h_Y = A^3 Y_1 + x_3 \overrightarrow{b} + x_2 A \overrightarrow{b} + x_1 A^2 \overrightarrow{e} + (A^2 + A + I) h_Y, \quad (20)$$

we have

$$(1 - A^3) Y_1 = x_3 \overrightarrow{b} + x_2 A \overrightarrow{b} + x_1 A^2 \overrightarrow{e} + (A^2 + A + I) h_y. \quad (21)$$

Therefore, if $A$ has no eigenvalue equal to 1, then we can find $Y_1$ explicitly as

$$Y_1 = (1 - A^3)^{-1} \left[ x_3 \overrightarrow{b} + x_2 A \overrightarrow{b} + x_1 A^2 \overrightarrow{e} + (A^2 + A + I) h_y \right]$$
$$= (1 - A^3)^{-1} \left[ \frac{1+a+ad}{1-a^2d} \hat{\mu} \overrightarrow{b} + \frac{1+d+ad}{1-a^2d} \hat{\mu} A \overrightarrow{b} + \frac{1+a+a^2}{1-a^2d} \hat{\mu} A^2 \overrightarrow{e} + (A^2 + A + I) h_y \right], \quad (22)$$

and periodic points $Z_1 = (x_1, Y_1)$, $Z_2 = (x_2, Y_2)$ and $Z_3 = (x_3, Y_3)$ of the 3-cycle $RL^2$ can be obtained explicitly. Also, for $\hat{\mu} > 0$ the parameter region for the existence of the 3-cycle $RL^2$, for which $x_1 > 0$ and $x_2, x_3 < 0$, is

$$\mathcal{R}_{\hat{\mu}>0} = \left\{ (a, d) \in \mathbb{R}^2 \,|\, a > 0, d < \frac{-a-1}{a} \right\}. \quad (23)$$

Now, for $\hat{\mu} > 0$ and $a, d \in \mathcal{R}_{\hat{\mu}>0}$, we show that the 3-cycle $RL^2$ appears through a border collision with symbolic sequence $RL0$ (where 0 denotes the point for which $x = 0$). For this purpose, consider a 3-cycle $RL^2$. Then, according to [25], the border collision arises if the condition $x_1 \leq \hat{\mu}$ holds. By (16), this condition results in $\frac{1+a+a^2}{1-a^2d} \leq 1$ as $\hat{\mu} > 0$. Then, $a > 0$ implies that $1 + a + ad \leq 0$, and the border collision occurs when $x_1 = \hat{\mu}$, i.e. $1 + a + ad = 0$. This means the 3-cycle $RL^2$ collides with the border for $x_3 = 0$ i.e. with symbolic sequence $RL0$. Moreover, as shown, the border collision bifurcation is independent of $\hat{\mu}$, see also [25]. In addition, the border collision bifurcation curve related to the 3-cycle $RL^2$ for $\hat{\mu} > 0$ is

$$\mathcal{C}_{\hat{\mu}>0} = \left\{ (a,d) \in \mathbb{R}^2 \,\big|\, a > 0,\, 1 + a + ad = 0 \right\}, \tag{24}$$

which is in fact the boundary of $\mathcal{R}_{\hat{\mu}>0}$. For $\hat{\mu} < 0$ similar results can be obtained as both cases, $\hat{\mu} > 0$ and $\hat{\mu} < 0$, are topologically conjugate. Indeed, for $\hat{\mu} < 0$ there is the parametric region $\mathcal{R}_{\hat{\mu}<0}$ with the boundary $\mathcal{C}_{\hat{\mu}<0}$ as its border collision bifurcation curve. Also, both parametric regions $\mathcal{R}_{\hat{\mu}>0}$ and $\mathcal{R}_{\hat{\mu}<0}$ are symmetric with respect to the diagonal $a = d$ (it is enough to change $(a,d)$ to $(d,a)$). The parametric regions $\mathcal{R}_{\hat{\mu}>0}$ and $\mathcal{R}_{\hat{\mu}<0}$ together with the bifurcation curves $\mathcal{C}_{\hat{\mu}>0}$ and $\mathcal{C}_{\hat{\mu}<0}$ are illustrated in Fig. 1.

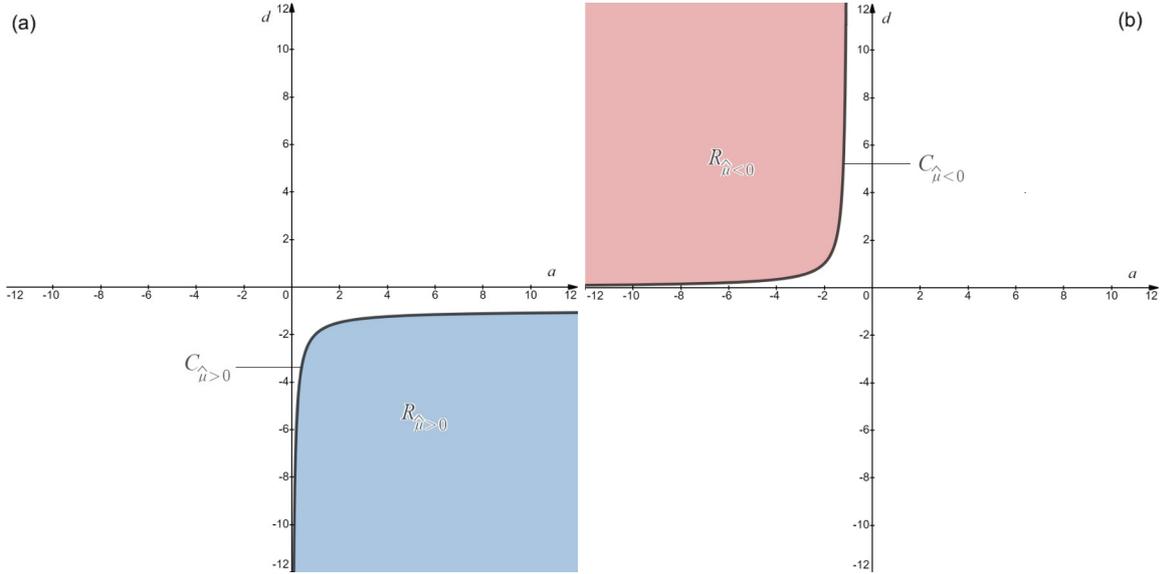

**Fig. 1** Parametric regions for the existence of 3-cycles $RL^2$ and $RLR$: (a) for $\hat{\mu} > 0$; (b) for $\hat{\mu} < 0$. The regions $\mathcal{R}_{\hat{\mu}>0}$ and $\mathcal{R}_{\hat{\mu}<0}$ are bounded by bifurcation curves $\mathcal{C}_{\hat{\mu}>0}$ and $\mathcal{C}_{\hat{\mu}<0}$, respectively.

In a likewise manner, we can also study the existence, stability and border collision bifurcations of all n-cycles $RL^{n-1}$ for any $n > 3$ for system (16)-(17). The following proposition addresses this issue.

**Proposition 1** *The n-cycle $\mathcal{O}_{RL^{n-1}}, n \geq 3$, of system (16)-(17) with $\hat{\mu} > 0$, exists if $A$ does not have any eigenvalue equal to 1 and $(a, d)$ belongs to the region*

$$\mathcal{R}_{\hat{\mu}>0}^{RL^{n-1}} = \left\{ (a,d) \in \mathbb{R}^2 \,\big|\, a > 0,\, d < -\frac{1-a^{n-1}}{(1-a)a^{n-2}} \right\}. \tag{25}$$

*Furthermore, a border collision bifurcation, independent of $\hat{\mu}$, with symbolic sequence $RL^{n-2}0$ occurs if $(a, d)$ lies on the bifurcation curve*

$$\mathcal{C}_{\hat{\mu}>0}^{RL^{n-1}} = \left\{ (a,d) \in \mathbb{R}^2 \,\big|\, a > 0,\, d = -\frac{1-a^{n-1}}{(1-a)a^{n-2}} \right\}, \tag{26}$$

such that $\mathcal{R}_{\hat{\mu}>0}^{RL^{n-1}}$ is bounded by the bifurcation curve $\mathcal{C}_{\hat{\mu}>0}^{RL^{n-1}}$. Similar results can be derived for $\hat{\mu} < 0$.

*Proof* Consider a basic n-cycle $\mathcal{O}_{RL^{n-1}}, n \geq 3$, of system (16)-(17) with periodic points $\{Z_1, Z_2, \ldots, Z_n\}$. Similar as in the previous section, $x$-components of the periodic orbits can be obtained from

$$\underbrace{f^- \circ \ldots \circ f^-}_{\text{n-1 times}} \circ f^+(x_1) = x_1 \tag{27}$$

as

$$x_1 = \frac{1-a^n}{(1-a)(1-a^{n-1}d)} \hat{\mu} > 0, \qquad x_2 = \frac{1-a+d(1-a^{n-1})}{(1-a)(1-a^{n-1}d)} \hat{\mu} < 0,$$

$$x_i = \frac{1-a^{i-1}+a^{i-2}d(1-a^{n-i+1})}{(1-a)(1-a^{n-1}d)} \hat{\mu} < 0, \qquad i = 3, 4, \cdots, n. \tag{28}$$

Also, due to

$$Y_2 = \vec{e} x_1 + AY_1 + h_Y,$$

$$Y_3 = \vec{b} x_2 + AY_2 + h_Y = \vec{b} x_2 + x_1 A \vec{e} + A^2 Y_1 + Ah_Y + h_Y,$$

$$Y_4 = \vec{b} x_3 + AY_3 + h_Y = A^3 Y_1 + x_3 \vec{b} + x_2 A \vec{b} + x_1 A^2 \vec{e} + (A^2 + A + I) h_Y,$$

$$Y_5 = \vec{b} x_4 + AY_4 + h_Y = A^4 Y_1 + x_4 \vec{b} + x_3 A \vec{b} + x_2 A^2 \vec{b} + x_1 A^3 \vec{e} + (A^3 + A^2 + A + I) h_Y,$$

$$\vdots$$

$$Y_1 = \vec{b} x_n + AY_n + h_Y = A^n Y_1 + x_n \vec{b} + x_{n-1} A \vec{b} + x_{n-2} A^2 \vec{b} + \ldots + x_2 A^{n-2} \vec{b}$$
$$+ x_1 A^{n-1} \vec{e} + (A^{n-1} + \ldots + A^2 + A + I) h_Y, \tag{29}$$

it is concluded that

$$(1 - A^n) Y_1 = x_n \vec{b} + x_{n-1} A \vec{b} + x_{n-2} A^2 \vec{b} + \ldots + x_2 A^{n-2} \vec{b} + x_1 A^{n-1} \vec{e}$$
$$+ (A^{n-1} + \ldots + A^2 + A + I) h_Y. \tag{30}$$

If $A$ has no eigenvalue equal to 1, then

$$Y_1 = (1 - A^n)^{-1} \big[ x_n \vec{b} + x_{n-1} A \vec{b} + x_{n-2} A^2 \vec{b} + \ldots + x_2 A^{n-2} \vec{b} + x_1 A^{n-1} \vec{e}$$
$$+ (A^{n-1} + \ldots + A^2 + A + I) h_Y \big]$$
$$= (1 - A^n)^{-1} \bigg[ \frac{1-a^{n-1}+a^{n-2}d(1-a)}{(1-a)(1-a^{n-1}d)} \hat{\mu} \vec{b} + \frac{1-a^{n-2}+a^{n-3}d(1-a^2)}{(1-a)(1-a^{n-1}d)} \hat{\mu} A \vec{b}$$
$$+ \frac{1-a^{n-3}+a^{n-4}d(1-a^3)}{(1-a)(1-a^{n-1}d)} \hat{\mu} A^2 \vec{b} + \cdots + \frac{1-a+d(1-a^{n-1})}{(1-a)(1-a^{n-1}d)} \hat{\mu} A^{n-2} \vec{b}$$
$$+ \frac{1-a^n}{(1-a)(1-a^{n-1}d)} \hat{\mu} A^{n-1} \vec{e} + (A^{n-1} + \ldots + A^2 + A + I) h_Y \bigg]. \tag{31}$$

Hence, if $A$ has no eigenvalue equal to 1, according to [25], for $\hat{\mu} > 0$ the cycle $\mathcal{O}_{RL^{n-1}}$ exists for all $(a, d)$ belonging to the region $\mathcal{R}_{\hat{\mu}>0}^{RL^{n-1}}$ given by (25). In this case, due to [25], the border collision bifurcation with symbolic sequence $RL^{n-1}0$, related to the collision of the periodic point $x_n = 0$, happens at $x_1 = \hat{\mu}$, provided that $(a, d)$ lies on the bifurcation curve $\mathcal{C}_{\hat{\mu}>0}^{RL^{n-1}}$ defined by (26). It is easy to see that $\mathcal{R}_{\hat{\mu}>0}^{RL^{n-1}}$ is bounded by the bifurcation curve $\mathcal{C}_{\hat{\mu}>0}^{RL^{n-1}}$.

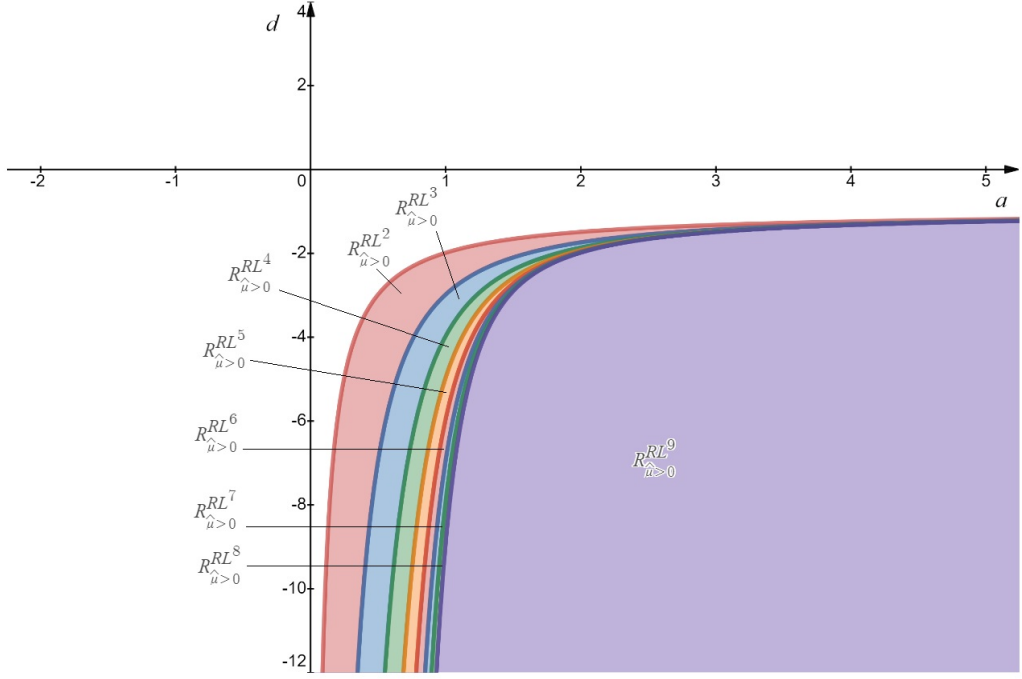

**Fig. 2** The parametric regions $\mathcal{R}_{\hat{\mu}>0}^{RL^{n-1}}$ for the existence of n-cycles $RL^{n-1}$, $n = 2, 3, \ldots, 9$, and the related bifurcation curves.

*Remark 1* As demonstrated in Fig. 2, for $n \geq 3$ the parameter regions $\mathcal{R}_{\hat{\mu}>0}^{RL^{n-1}}$, for which the n-cycles $RL^{n-1}$ exist, satisfy

$$\mathcal{R}_{\hat{\mu}>0}^{RL^2} \supset \mathcal{R}_{\hat{\mu}>0}^{RL^3} \supset \mathcal{R}_{\hat{\mu}>0}^{RL^4} \supset \ldots \supset \mathcal{R}_{\hat{\mu}>0}^{RL^i} \supset \ldots \qquad (32)$$

Note that for $\hat{\mu} < 0$ there are similar relations between the parameter regions $\mathcal{R}_{\hat{\mu}<0}^{RL^{n-1}}, n \geq 3$.

**Corollary 1** *If all eigenvalues of A are smaller than 1 in modulus, then the stability of the n-cycles $\mathcal{O}_{RL^{n-1}}$ is determined by the map (16). This means that the n-cycle $\mathcal{O}_{RL^{n-1}}$ is attracting if $(a,d)$ belongs to the parametric region*

$$\mathcal{R}_s^{RL^{n-1}} = \left\{ (a,d) \in \mathbb{R}^2 \,\Big|\, a > 0, \; -\frac{1}{a^{n-1}} < d < -\frac{1-a^{n-1}}{(1-a)a^{n-2}} \right\}. \qquad (33)$$

*Proof* Since $(x_1, Y_1)$ is the fixed point of the map

$$\underbrace{F^- \circ \ldots \circ F^-}_{\text{n-1 times}} \circ F^+(x,Y) = \begin{pmatrix} a^{n-1} d\, x \\ A^n Y \end{pmatrix} \qquad (34)$$

$$+ \begin{pmatrix} \dfrac{1-a^n}{1-a}\hat{\mu} \\ x_n \overrightarrow{b} + x_{n-1} A \overrightarrow{b} + x_{n-2} A^2 \overrightarrow{b} + \ldots + x_2 A^{n-2} \overrightarrow{b} + x_1 A^{n-1} \overrightarrow{e} + (A^{n-1} + \ldots + A^2 + A + I)\, h_Y \end{pmatrix}, \qquad (35)$$

so $\mathcal{O}_{RL^{n-1}}$ is an attracting cycle, provided that $A$ has only eigenvalues smaller in modulus than 1 and $-1 < a^{n-1}d$ ($(a,d) \in \mathcal{R}_{\hat{\mu}>0}^{RL^{n-1}}$, hence $a^{n-1}d < 0$). Due to (25), the last assertion follows if $(a,d)$ belongs to the parametric region $\mathcal{R}_s^{RL^{n-1}}$. Please also see [1,25] for more details on the parametric regions for stability of cycles.

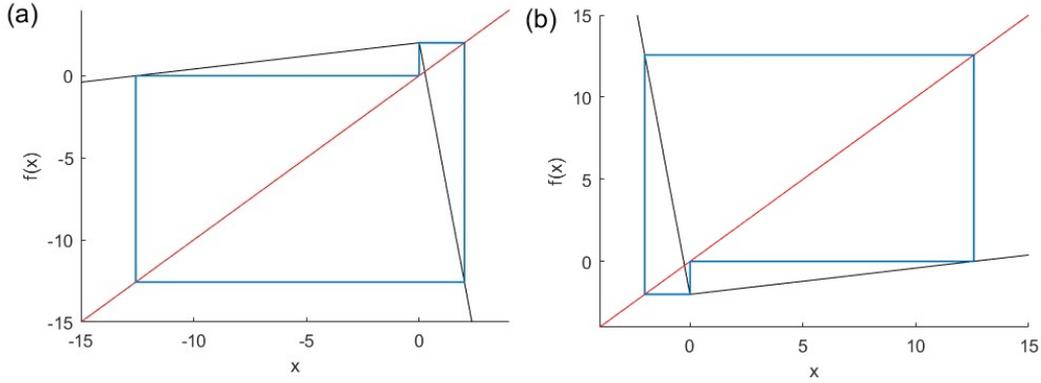

**Fig. 3** (a) The stable 3-cycle $RL^2$ for $\hat{\mu} = 2$ and $(a,d) = (0.16, -7.29)$ with the periodic points $\{1.9982, -12.5674, -0.0107\}$; (b) The stable 3-cycle $RLR$ for $\hat{\mu} = -2$ and $(a,d) = (-7.29, 0.16)$ with the periodic points $\{0.0107, -1.9982, 12.5674\}$.

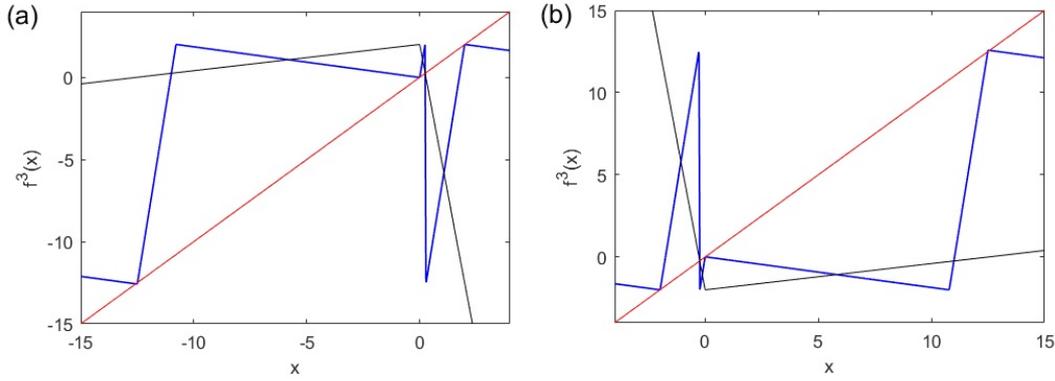

**Fig. 4** The skew tent map (16) and its third iteration for: (a) the 3-cycle $RL^2$; and (b) the 3-cycle $RLR$ of Fig. 3.

*Remark 2* For continuous map (16)-(17), any basic n-cycle $RL^{n-1}$, $n \geq 3$, exists in a pair with its complementary n-cycle with symbolic sequence $RL^{n-2}R$. Also, the border collision of the n-cycles $RL^{n-2}R$ occurs analogously with the same bifurcation curve (26), see [1, 25]. Moreover, as shown in corollary 1, for $\hat{\mu} > 0$, there is a stability region in which only the cycle $RL^{n-1}$ can be stable, whereas its complementary cycle is always unstable. Note that for $\hat{\mu} < 0$, there exists a stability region within which only the cycle $RL^{n-2}R$ can be stable, as in this case $RL^{n-1}$ is always unstable. For example, for $\hat{\mu} = 2$ and $(a,d) = (0.16, -7.29) \in \mathcal{R}_s^{RL^2}$, there is the stable 3-cycle $RL^2$ for (16)-(17) shown in Fig. 3(a). On the other hand, for $\hat{\mu} = -2$ and the symmetric values $(a,d) = (-7.29, 0.16)$, the attracting cycle is $RLR$ (Fig. 3(b)). Furthermore, these two cycles have the same periodic points, but with different signs; see Fig. 3 and Fig. 4.

*Remark 3* By theorem 1, the existence of a period-3 orbit for one-dimensional continuous maps implies the existence of chaos. The particular system considered here is indeed a one-dimensional continuous system in variable $x$, coupled to a $(n-1)$-dimensional system of similar form, such that the conditions of theorem 1 apply. In particular, when all eigenvalues of $A$ are $\neq 1$, according to [1, 25], there are the cyclic $n$-band and $2n$-band chaotic attractors $\mathcal{A}_{n,n}$ and $\mathcal{A}_{n,2n}$ for all $(a,d)$ belonging to $\mathcal{R}_{n,n}^{\mathcal{A}_{n,n}}$ and $\mathcal{R}_{n,2n}^{\mathcal{A}_{n,2n}}$, respectively, where

$$\mathcal{R}_{n,n}^{\mathcal{A}_{n,n}} = \left\{ (a,d) \in \mathbb{R}^2 \mid d < -\frac{1-a^{n-1}}{(1-a)a^{n-2}}, \ a^{2(n-1)}d^3 + a - d < 0, \ a^{n-1}d^2 + d - a < 0 \right\}, \quad (36)$$

$$\mathcal{R}_{n,2n}^{\mathcal{A}_{n,2n}} = \left\{ (a,d) \in \mathbb{R}^2 \mid d < -\frac{1-a^{n-1}}{(1-a)a^{n-2}}, \ d < -\frac{1}{a^{n-1}}, \ a^{2(n-1)}d^3 + a - d > 0 \right\}. \quad (37)$$

## 4 Application of results to PLRNNs

In this section we discuss how the results obtained in the previous section can be extended to PLRNNs.

Consider the discrete-time PLRNN

$$Z_{t+1} = (A + WD_{\Omega(t)})Z_t + h, \tag{38}$$

where $Z_t = (z_{1t}, \cdots, z_{Mt})^T \in \mathbb{R}^M$ represents the neural state vector at time $t = 1...T$, $A = diag(a_{11}, \cdots, a_{MM}) \in \mathbb{R}^{M \times M}$ is a diagonal matrix of auto-regression weights, $W \in \mathbb{R}^{M \times M}$ denotes an off-diagonal matrix of connection weights with diagonal elements equal to zero, and $h$ is a bias term [6,14]. Also,

$$D_{\Omega(t)} := \text{diag}(d_{\Omega(t)}),$$

with

$$d_{\Omega(t)} := \big(d_1(t), d_2(t), \cdots, d_M(t)\big),$$

such that $d_i(t) = 0$ if $z_{it} \leq 0$ and $d_i(t) = 1$ if $z_{it} > 0$, for $i = 1, 2, \cdots, M$. Based on the sign of the $M$ components of $Z_t$, there are $2^M$ distinct configurations for matrix $D_{\Omega(t)}$. That means, the phase space of system (38) is divided into $2^M$ sub-regions by $M2^{M-1}$ hyper-surfaces as switching borders. If we index the $2^M$ different configurations of $D_{\Omega(t)}$ as $D_{\Omega^k}$, $k \in \{1, 2, \cdots, 2^M\}$, we can define $2^M$ matrices

$$W_{\Omega^k} := A + WD_{\Omega^k}, \tag{39}$$

such that in each sub-region $S_{\Omega^k}$, $k \in \{1, 2, \cdots, 2^M\}$, the dynamics are described by a different map

$$Z_{t+1} = F(Z_t) = W_{\Omega^k} Z_t + h, \qquad Z_t \in S_{\Omega^k}. \tag{40}$$

The sub-regions $S_{\Omega^k}$'s can be defined as [20]

$$\begin{cases} S_{\Omega^1} = \hat{S}_0 = \hat{S}_{(\underbrace{0\,0\,0\,\cdots\,0}_{M})_2^*} = \hat{S}_{\underbrace{0\,0\,0\,\cdots\,0}_{M}} = \Big\{Z_t \in \mathbb{R}^M; z_{it} \leq 0,\, i=1,2,\cdots,M\Big\}, \\ S_{\Omega^2} = \hat{S}_1 = \hat{S}_{(\underbrace{0\,0\,\cdots\,0\,1}_{M})_2^*} = \hat{S}_{\underbrace{1\,0\,0\,\cdots\,0}_{M}} = \Big\{Z_t \in \mathbb{R}^M; z_{1t} > 0, z_{it} \leq 0,\, i \neq 1\Big\}, \\ S_{\Omega^3} = \hat{S}_2 = \hat{S}_{(\underbrace{0\,\cdots\,0\,1\,0}_{M})_2^*} = \hat{S}_{\underbrace{0\,1\,0\,\cdots\,0}_{M}} = \Big\{Z_t \in \mathbb{R}^M; z_{2t} > 0, z_{it} \leq 0,\, i \neq 2\Big\}, \\ S_{\Omega^4} = \hat{S}_3 = \hat{S}_{(\underbrace{0\,\cdots\,0\,1\,1}_{M})_2^*} = \hat{S}_{\underbrace{1\,1\,0\,\cdots\,0}_{M}} = \Big\{Z_t \in \mathbb{R}^M; z_{1t}, z_{2t} > 0, z_{it} \leq 0,\, i \neq 1,2\Big\}, \\ S_{\Omega^5} = \hat{S}_4 = \hat{S}_{(\underbrace{0\,\cdots\,1\,0\,0}_{M})_2^*} = \hat{S}_{\underbrace{0\,0\,1\,0\,\cdots\,0}_{M}} = \Big\{Z_t \in \mathbb{R}^M; z_{3t} > 0, z_{it} \leq 0,\, i \neq 3\Big\}, \\ \vdots \qquad\qquad\qquad \vdots \\ S_{\Omega^{2^M}} = \hat{S}_{2^M-1} = \hat{S}_{(\underbrace{1\,1\,1\,\cdots\,1}_{M})_2^*} = \hat{S}_{\underbrace{1\,1\,1\,\cdots\,1}_{M}} = \Big\{Z_t \in \mathbb{R}^M; z_{it} > 0,\, i=1,2,\cdots,M\Big\}. \end{cases} \tag{41}$$

where each subindex $d$ of $\hat{S}$, $0 \leq d \leq 2^M - 1$, is associated with a sequence $d_M\, d_{M-1} \cdots d_2\, d_1$ of binary digits. The notation $(d_1\, d_2 \cdots d_M)_2^*$ in building each corresponding sequence stands for the mirror image of the binary representation of $d$ with $M$ digits. By mirror image here we mean writing digits $d_1\, d_2 \cdots d_M$ from right to left, i.e. $d_M\, d_{M-1} \cdots d_2\, d_1$.

Denoting switching boundaries $\Sigma_{ij} = \bar{S}_{\Omega^i} \cap \bar{S}_{\Omega^j}$ between every pair of successive sub-regions $S_{\Omega^i}$ and $S_{\Omega^j}$ with $i, j \in \{1, 2, \cdots, 2^M\}$, we can rewrite map (38) as

$$Z_{t+1} = F(Z_t) = \begin{cases} F_1(Z_t) = W_{\Omega^1} Z_t + h; & Z_t \in \bar{S}_{\Omega^1} \\ F_2(Z_t) = W_{\Omega^2} Z_t + h; & Z_t \in \bar{S}_{\Omega^2} \\ F_3(Z_t) = W_{\Omega^3} Z_t + h; & Z_t \in \bar{S}_{\Omega^3} \\ F_4(Z_t) = W_{\Omega^4} Z_t + h; & Z_t \in \bar{S}_{\Omega^4} \\ \vdots & \vdots \\ F_{2^M}(Z_t) = W_{\Omega^{2^M}} Z_t + h; & Z_t \in \bar{S}_{\Omega^{2^M}} \end{cases}. \quad (42)$$

Hence, $S_{\Omega^i}$ and $S_{\Omega^j}$ are two successive sub-regions with the switching boundary $\Sigma_{ij} = \bar{S}_{\Omega^i} \cap \bar{S}_{\Omega^j}$, iff there is exactly one $1 \leq s \leq M$ such that for all $(z_{i_1 t}, \cdots, z_{i_M t})^T \in \mathring{S}_{\Omega^i}$ and $(z_{j_1 t}, \cdots, z_{j_M t})^T \in \mathring{S}_{\Omega^j}$

$$\begin{cases} z_{i_s t} \cdot z_{j_s t} < 0 \\ z_{i_r t} \cdot z_{j_r t} > 0, \; 1 \leq \underset{r \neq s}{r} \leq M \end{cases}, \quad (43)$$

which, for $i_k = j_k = k$, $k = 1, \cdots, M$, is equivalent to

$$\begin{cases} z_{st} \cdot z_{st} < 0 \\ z_{rt} \cdot z_{rt} > 0, \; 1 \leq \underset{r \neq s}{r} \leq M \end{cases}, \quad (44)$$

Now, for $s = 1$, we can extend the previous results for PLRNNs. For this purpose, we investigate n-cycles of system (42) locally near only one boundary $\Sigma_{ij}$. That means, restricting the domain of function $F$ on $U$, denoted by $G = F|_U$, we consider the system

$$Z_{t+1} = G(Z_t) = \begin{cases} G_i(Z_t) = W_{\Omega^i} Z_t + h; & Z_t \in \overline{S_{\Omega^i} \cap U} := \hat{S}_{\Omega^i}(z_{1t} \leq 0) \\ G_j(Z_t) = W_{\Omega^j} Z_t + h; & Z_t \in \overline{S_{\Omega^j} \cap U} := \hat{S}_{\Omega^j}(z_{1t} \geq 0) \end{cases}, \quad (45)$$

where $U$ is an open subset of $\mathbb{R}^M$ such that $U \cap \Sigma_{ij} \neq \varnothing$ and $0 \notin U \cap \Sigma_{ij}$. Now, assume $W_{\Omega^i} = A_1$ and $W_{\Omega^j} = A_2$ in system (45), where $A_1$ and $A_2$ are in the form (9). Then, writing $Z = (x, Y)$, where $x$ indicates the first component of $Z$ and $Y$ is a vector collecting the remaining $(M - 1)$ components, we can use the previous results to establish a necessary condition for the existence of the n-cycles $RL^{n-1}$, $n \geq 3$, and their associated border collision bifurcations, for (45). Also, we can extend the results obtained for system (45). More precisely, let $\hat{\mu} := h_1 > 0$ and $\mathcal{O}_{RL^{n-1}}$ be an n-cycle for system (45) with periodic fixed points $\{(x_1, Y_1), (x_2, Y_2), \cdots, (x_n, Y_n)\}$. If none of the eigenvalues of $A$ are equal to 1, then, as it was shown in the previous section, we can obtain $Y_1, Y_2, \cdots, Y_n$ explicitly. In this case, the n-cycles $RL^{n-1}$, $n \geq 3$, exist if $a, d \in \mathcal{R}_{\hat{\mu}>0}^{RL^{n-1}}$. Other results derived in the previous sections can be extended in a similar way.

The results presented here are an important step toward a more systematic study of cycles in PLRNNs as they have been found, for instance, in PLRNNs inferred from brain imaging data [14]. In this context we also remark that the conditions used in Lemma 1, in particular $\vec{c} = \vec{0}$, are not too restrictive, as recent results on the application of PLRNNs to diverse challenging machine learning benchmarks and data sets have shown that zeroing out the 'inputs' from other states for a subset of PLRNN variables (i.e., setting row vectors of the transition matrix to 0 except for the $s$-th entry) can lead to much improved performance [24]. Specifically, units for which the self-connection is close to 1 while those conveying inputs form other units are close to 0 essentially become memory units which integrate external inputs, and thus help the PLRNN to solve problems that contain long-range dependencies and very slow time scales. Our results may thus not only help to understand and analyze the cyclic behavior of trained PLRNNs, but also bifurcations in these systems as parameters change throughout the training process [5].

## 5 Numerical simulations

*Example 1* Suppose that for the system in (16)-(17), $A = diag(a_1, a_2, \ldots, a_{n-1})$ is a diagonal matrix such that $a_i \neq 1$ for all $i \in \{1, 2, \ldots, n-1\}$. By proposition 1, for $\hat{\mu} > 0$, there exists a k-cycle $\mathcal{O}_k$ with symbolic form $RL^{k-1}$, if $(a, d) \in \mathcal{R}_{\hat{\mu}>0}^{RL^{k-1}}$. In this case, the periodic points $\{(x_1, Y_1), (x_2, Y_2), \ldots, (x_k, Y_k)\}$ can be calculated according to the proof of proposition 1. Specifically, $x_1, x_2, \ldots, x_k$ and $Y_1, Y_2, \ldots, Y_k$ can be obtained by (28) and (29), respectively. Also, since $a_i \neq 1$ for all $i \in \{1, 2, \ldots, n-1\}$, $Y_1$ has the explicit form (31), and thereby $Y_2, \ldots, Y_k$ can also be computed explicitly. Moreover, let $\vec{b} = (b_1, b_2, \ldots, b_{n-1})^T$, $\vec{e} = (e_1, e_2, \ldots, e_{n-1})^T$, and $h_Y = (h_Y^1, h_Y^2, \ldots, h_Y^{n-1})^T$, then, by diagonality of $A$, it is possible to simplify $Y_1$ as

$Y_1 =$

$$\begin{pmatrix} \frac{1}{1-a_1^k}\left[\left(x_k + x_{k-1} a_1 + x_{k-2} a_1^2 + \ldots + x_2 a_1^{k-2}\right) b_1 + x_1 a_1^{k-1} e_1 + \left(1 + a_1 + a_1^2 + \ldots + a_1^{k-1}\right) h_Y^1\right] \\ \frac{1}{1-a_2^k}\left[\left(x_k + x_{k-1} a_2 + x_{k-2} a_2^2 + \ldots + x_2 a_2^{k-2}\right) b_2 + x_1 a_2^{k-1} e_2 + \left(1 + a_2 + a_2^2 + \ldots + a_2^{k-1}\right) h_Y^2\right] \\ \vdots \\ \frac{1}{1-a_{n-1}^k}\left[\left(x_k + x_{k-1} a_{n-1} + \ldots + x_2 a_{n-1}^{k-2}\right) b_{n-1} + x_1 a_{n-1}^{k-1} e_{n-1} + \left(1 + a_{n-1} + \ldots + a_{n-1}^{k-1}\right) h_Y^{n-1}\right] \end{pmatrix}.$$
(46)

Now for $n = 4$, i.e. the 4-dimensional system (16)-(17), we have $A = diag(a_1, a_2, a_3)$, $\vec{b} = (b_1, b_2, b_3)^T$, $\vec{e} = (e_1, e_2, e_3)^T$, $Y = (Y^1, Y^2, Y^3)^T$ and $h_Y = (h_Y^1, h_Y^2, h_Y^3)^T$. Then, for $\hat{\mu} > 0$ and $a_1, a_2, a_3 \neq 1$, by proposition 1, there exists a 3-cycle $\mathcal{O}_3$, with symbolic sequence $RL^2$, if $(a, d) \in \mathcal{R}_{\hat{\mu}>0}^{RL^2} = \mathcal{R}_{\hat{\mu}>0}$. Furthermore, $O_3 = \{(x_1, Y_1), (x_2, Y_2), (x_3, Y_3)\}$, such that $x_1, x_2, x_3$ are given by (18)-(19). Hence, using (29) and (46) we get

$$Y_1 = \begin{pmatrix} Y_1^1 \\ Y_1^2 \\ Y_1^3 \end{pmatrix} = \begin{pmatrix} \frac{1}{1-a_1^3}\left[(x_3 + x_2 a_1) b_1 + x_1 a_1^2 e_1 + (1 + a_1 + a_1^2) h_Y^1\right] \\ \frac{1}{1-a_2^3}\left[(x_3 + x_2 a_2) b_2 + x_1 a_2^2 e_2 + (1 + a_2 + a_2^2) h_Y^2\right] \\ \frac{1}{1-a_3^3}\left[(x_3 + x_2 a_3) b_3 + x_1 a_3^2 e_3 + (1 + a_3 + a_3^2) h_Y^3\right] \end{pmatrix},$$

$$Y_2 = \begin{pmatrix} e_1 x_1 + a_1 Y_1^1 + h_Y^1 \\ e_2 x_1 + a_2 Y_1^2 + h_Y^2 \\ e_3 x_1 + a_3 Y_1^3 + h_Y^3 \end{pmatrix}, \qquad Y_3 = \begin{pmatrix} b_1 x_2 + x_1 a_1 e_1 + a_1^2 Y_1^1 + (1+a_1)h_Y^1 \\ b_2 x_2 + x_1 a_2 e_2 + a_2^2 Y_1^2 + (1+a_2)h_Y^2 \\ b_3 x_2 + x_1 a_3 e_3 + a_3^2 Y_1^3 + (1+a_3)h_Y^3 \end{pmatrix}. \qquad (47)$$

Let us choose the parameters of the system as follows:

$$a = 0.4, \quad d = -4, \quad \vec{b} = (1, 0.5, 0.6)^T \quad \vec{e} = (0.5, 1, 1)^T, \quad A = diag(0.4, 0.5, 0.6),$$
$$h_Y = (1, 0, 1)^T, \quad \hat{\mu} = 0.8. \tag{48}$$

Since $a_1, a_2, a_3 \neq 1$ and $(a, d) \in \mathcal{R}_{\hat{\mu}>0}$, there is a 3-cycle denoted by $\mathcal{O}_3 = \{(x_1, Y_1)^T, (x_2, Y_2)^T, (x_3, Y_3)^T\}$, where

$$(x_1, Y_1)^T = (x_1, Y_1^1, Y_1^2, Y_1^3)^T = (0.7610, 0.6685, -0.4794, 1.7444)^T,$$
$$(x_2, Y_2)^T = (x_2, Y_2^1, Y_2^2, Y_2^3)^T = (-2.2439, 1.6479, 0.5213, 2.8076)^T,$$
$$(x_3, Y_3)^T = (x_3, Y_3^1, Y_3^2, Y_3^3)^T = (-0.0976, -0.5847, -0.8613, 1.3382)^T.$$

In addition, the eigenvalues of $\mathcal{O}_3$ are $\lambda_1 = a^2 d = -0.64$, $\lambda_2 = a_1 = 0.4$, $\lambda_3 = a_2 = 0.5$, $\lambda_4 = a_3 = 0.6$, which are all less than one in magnitude. Also, $(a, d) \in \mathcal{R}_s^{RL^2}$ which confirms the stability of $\mathcal{O}_3$ due to corollary 1. Therefore, $\mathcal{O}_3$ is a stable 3-cycle for system (16)-(17), as illustrated in Fig. 5.

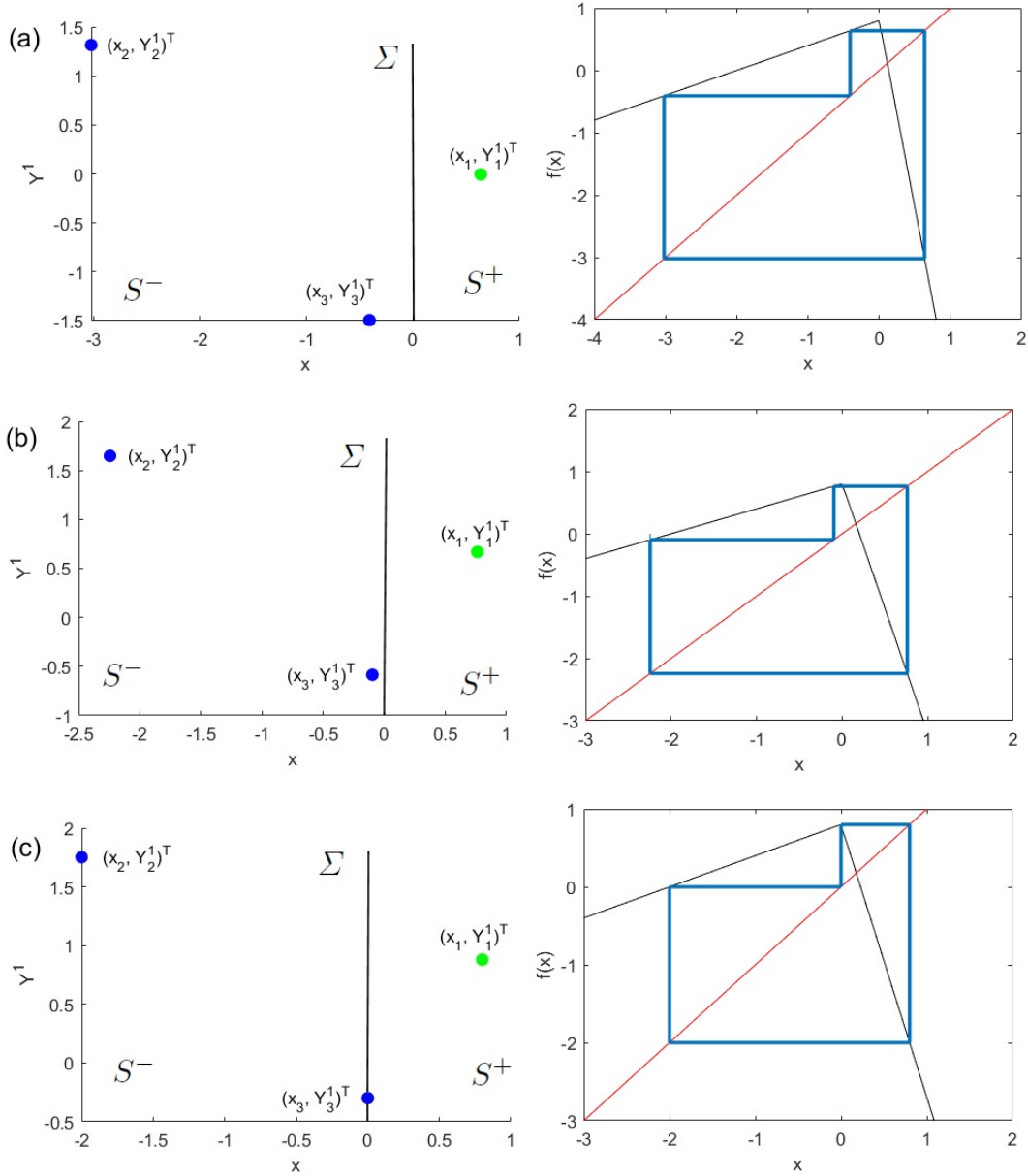

**Fig. 5** Example of a stable 3-cycle $\mathcal{O}_3$ of system (16)-(17) for $(a,d) \in \mathcal{R}_{\hat{\mu}>0}$. (a) $d = -6$, $a = 0.4$. (b) $d = -4$, $a = 0.4$. (c) Occurrence of a border-collision bifurcation with the periodic point $x_3 = 0$ for $(a,d) = (0.4, -3.5) \in \mathcal{C}_{\hat{\mu}>0}$.

Now, let the parameter $a$ be fixed to $a = 0.4$ and $d$ free to vary. As long as $d < \frac{-a-1}{a} = -3.5$, i.e. $(a,d) \in \mathcal{R}_{\hat{\mu}>0}$, there will be a 3-cycle for the system such that its periodic point $(x_3, Y_3)$ moves toward the border as $d$ tends to $-3.5$. Also at $d = d_{bif} = -3.5$, we have $1 + a + ad = 0$, i.e. $(a,d) \in \mathcal{C}_{\hat{\mu}>0}^{RL^2} = \mathcal{C}_{\hat{\mu}>0}$, from which by proposition 1 we can deduce the existence of a border-collision bifurcation with symbolic sequence $RL0$. Denoting this 3-cycle $RL0$ by $\mathcal{O}_3^c = \{(x_1, Y_1)^T, (x_2, Y_2)^T, (x_3, Y_3)^T\}$, we

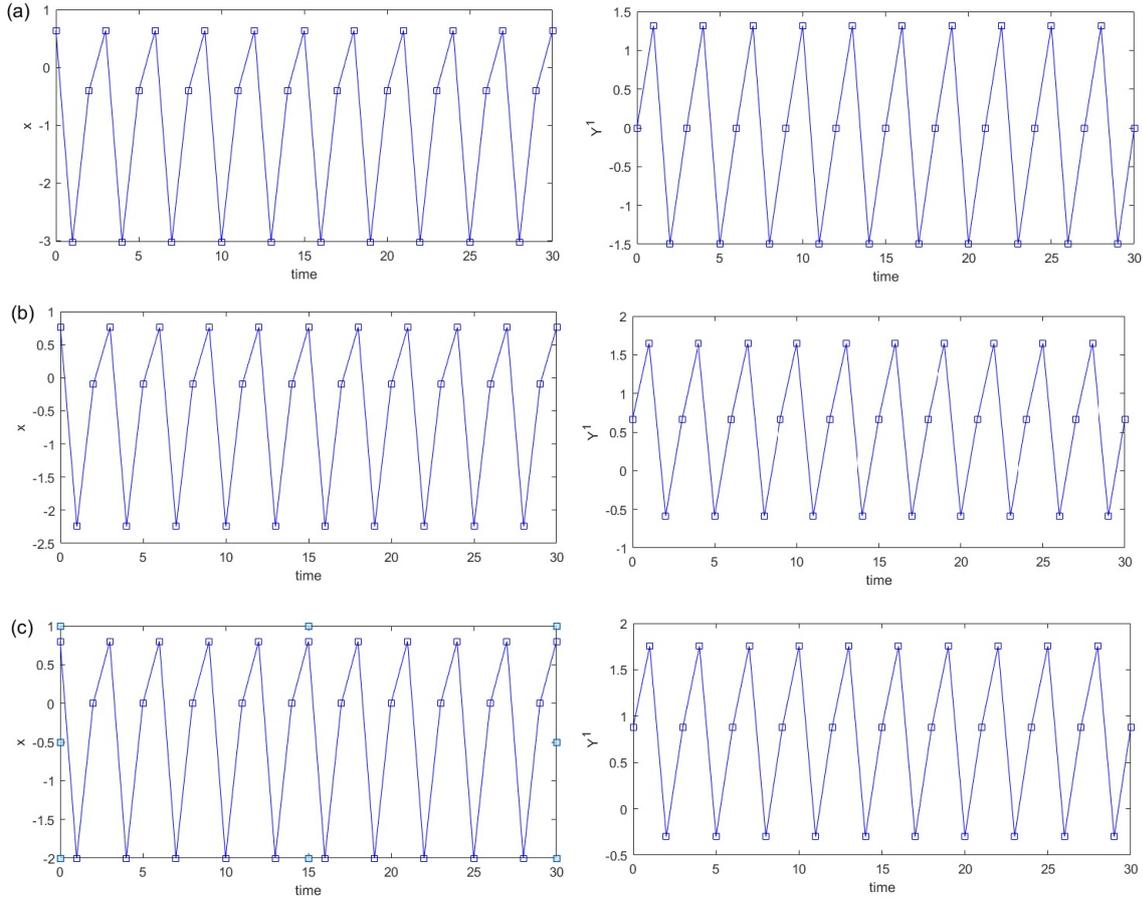

**Fig. 6** Time history of the stable 3-cycle $\mathcal{O}_3$ for $(a,d) \in \mathcal{R}_{\hat{\mu}>0}$; (a) $d = -6$, $a = 0.4$; (b) $d = -4$, $a = 0.4$; (c) $(a,d) = (0.4, -3.5) \in \mathcal{C}_{\hat{\mu}>0}$.

find

$$(x_1, Y_1)^T = (0.8,\ 0.8803,\ -0.3429,\ 1.9490)^T,$$
$$(x_2, Y_2)^T = (-2,\ 1.7521,\ 0.6286,\ 2.9694)^T,$$
$$(x_3, Y_3)^T = (0,\ -0.2991,\ -0.6857,\ 1.5816)^T.$$

In this case, a border-collision bifurcation happens at $x_1^c = \hat{\mu} = 0.8$ with the periodic point $x_3^c = 0$ (Fig. 5). Fig. 6 furthermore shows the time graph of the 3-cycle $\mathcal{O}_3$ for different values of $d$. For $d > d_{bif}$, on the other hand, there is no 3-cycle, as in this case $a, d \notin \mathcal{R}_{\hat{\mu}>0}$, (Fig. 7). Fig. 8 illustrates cobweb diagrams for these orbits. Finally, due to remark 3, for $(a, d) = (0.4, -6.5) \in \mathcal{R}_{3,3}^{\mathcal{A}_{3,3}}$ there is the cyclic 3-band chaotic attractor $\mathcal{A}_{3,3}$, and for $(a, d) = (0.4, -6.4) \in \mathcal{R}_{3,6}^{\mathcal{A}_{3,6}}$ the cyclic 6-band chaotic attractor $\mathcal{A}_{3,6}$ (Fig. 9). Finally, other n-cycles, with $n > 3$, can be obtained using the results in subsection 3.1, see Fig. 10.

## 6 Conclusions

Our aim in this paper was to investigate n-cycles $RL^{n-1}(n \geq 3)$ and border-collision bifurcations of a class of the piecewise linear continuous map (5) on $\mathbb{R}^n$. First, we considered a special partitioned form for the matrices $A_1$ and $A_2$ of the map (5). Under these conditions, we could determine some

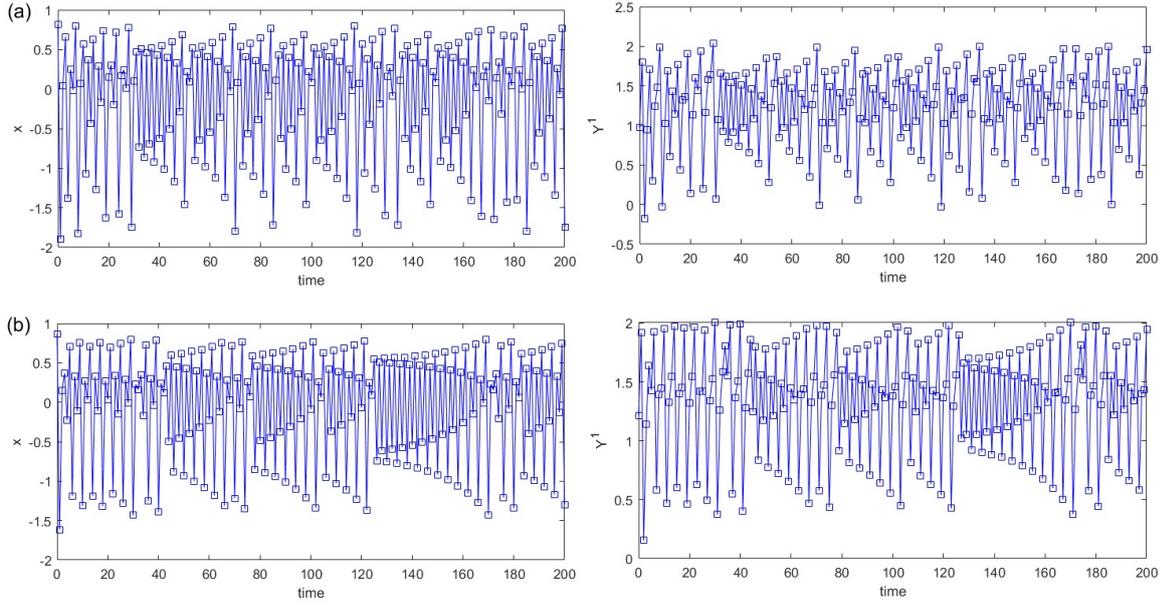

**Fig. 7** Time histories for two orbits of system (16)-(17) for $(a,d) \notin \mathcal{R}_{\hat{\mu}>0}$ and $d > d_{bif}$: (a) $d = -3.3$, $a = 0.4$; (b) $d = -2.8$, $a = 0.4$.

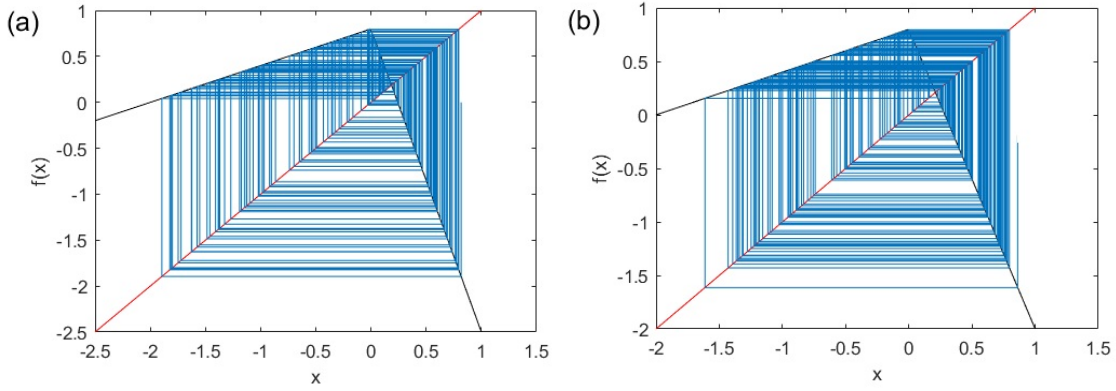

**Fig. 8** Cobweb plot of system (16)-(17) for $(a,d) \notin \mathcal{R}_{\hat{\mu}>0}$ and $d > d_{bif}$: (a) $d = -3.3$, $a = 0.4$; (b) $d = -2.8$, $a = 0.4$.

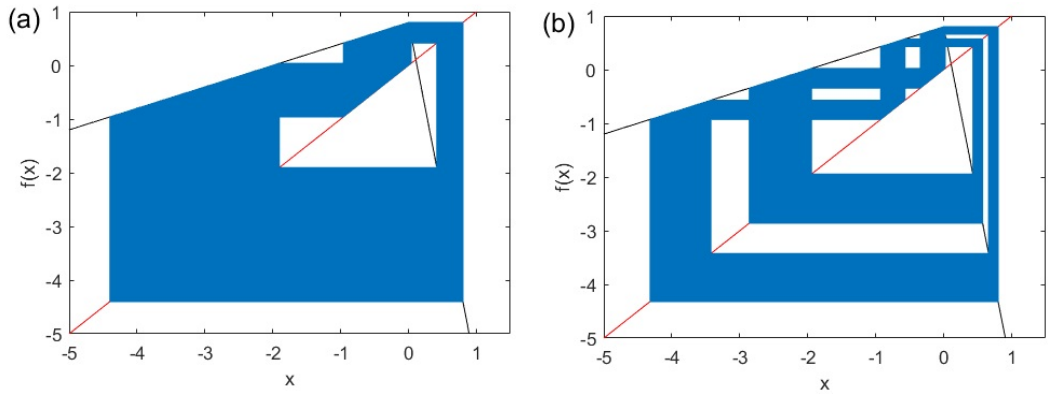

**Fig. 9** Cyclic chaotic attractors of system (16)-(17): (a) the cyclic 3-band chaotic attractor $\mathcal{A}_{3,3}$ for $(a,d) = (0.4, -6.5) \in \mathcal{R}_{3,3}^{\mathcal{A}_{3,3}}$; (b) the cyclic 6-band chaotic attractor $\mathcal{A}_{3,6}$ for $(a,d) = (0.4, -6.4) \in \mathcal{R}_{3,6}^{\mathcal{A}_{3,6}}$.

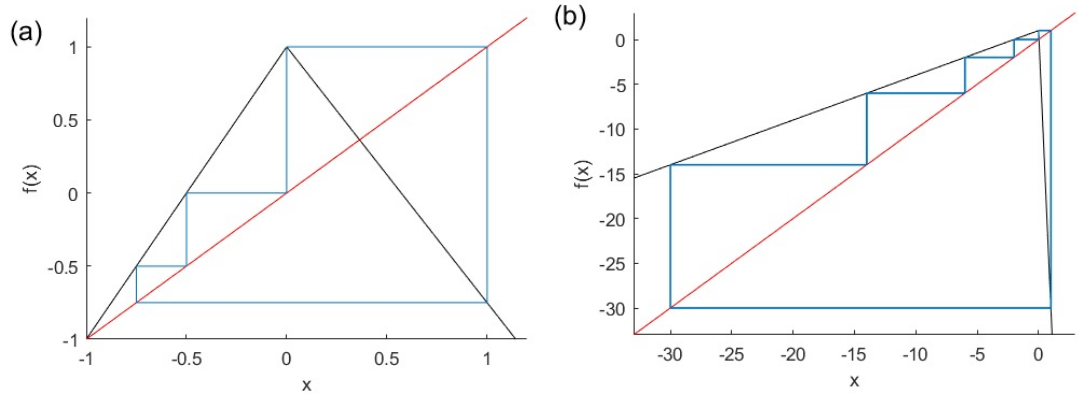

**Fig. 10** (a) The unstable 4-cycle $RL^20$ for $\hat{\mu} = 1$ and $(a, d) = (2, -1.75) \in \mathcal{C}_{\hat{\mu}>0}^{RL^3}$ corresponding to the border collision bifurcation; (b) The stable 6-cycle $RL^5$ for $\hat{\mu} = 1$ and $(a, d) = (0.5, -31) \in \mathcal{R}_s^{RL^5}$.

parametric regions for the existence (or nonexistence) of 3-cycles and occurrence of border collision bifurcations. Furthermore, we obtained some conditions for the existence of n-cycles $RL^{n-1}(n > 3)$ and their border collision bifurcations in system (16)-(17), and demonstrated some of the results by numerical simulations. Finally, we indicated how our findings can be applied to PLRNNs.

The results presented here are an important step toward a more systematic study of cycles in PLRNNs. Such cycles have been demonstrated, for instance, in PLRNNs inferred from brain imaging or other types of experimental data [14, 24]. In this context we also reiterate, as stated above, that the conditions used in Lemma 1, in particular $\vec{c} = \vec{0}$, do not curtail the application of our results to PLRNNs too much, as similar conditions may in fact enhance the performance of PLRNNs, and facilitate PLRNN training, on a variety of machine learning and dynamical systems identification problems [24]. Our results may not only help to understand and analyze the cyclic behavior of trained PLRNNs, but also bifurcations in these systems as parameters change throughout the training process [5].

## 7 Acknowledgments

This work was funded by a grant from the German Science Foundation (DFG) to DD (Du 354/10-1). We are very grateful to one anonymous referee who provided specific ideas for rewriting our system, and pointed out to us how our results could be improved and generalized.

## 8 Conflict of Interest

The authors declare that they have no conflict of interest.